\documentclass[graybox]{svmult}

\usepackage{amsfonts,amsmath,amstext,latexsym,amssymb}
\usepackage{mathptmx}       
\usepackage{helvet}         
\usepackage{courier}        
\usepackage{type1cm}        
\usepackage{makeidx}         
\usepackage{graphicx}        
\usepackage{multicol}        
\usepackage[bottom]{footmisc}

\makeindex             

\usepackage{cite}

\usepackage[inline]{enumitem}  
\setlist{labelindent=1pt,itemsep=.5em}
\setlist[itemize]{leftmargin=1cm}
\setlist[enumerate]{itemindent=0em,leftmargin=1cm}

\usepackage[draft=false,hidelinks,colorlinks=true,linkcolor=blue,citecolor=black,urlcolor=black,anchorcolor=blue,bookmarks=false,hypertexnames=true,backref=page]{hyperref}

\allowdisplaybreaks[3]

\smartqed

\newtheorem{rem:DjinjaTumwSilv1}{Remark}[section]
\newtheorem{lem:DjinjaTumwSilv1}{Lemma}[section]
\newtheorem{thm:DjinjaTumwSilv1}{Theorem}[section]
\newtheorem{cor:DjinjaTumwSilv1}[thm:DjinjaTumwSilv1]{Corollary}
\newtheorem{prop:DjinjaTumwSilv1}[thm:DjinjaTumwSilv1]{Proposition}
\newtheorem{ex:DjinjaTumwSilv1}[thm:DjinjaTumwSilv1]{Example}

\DeclareMathOperator{\esssup}{ess\;sup}
\DeclareMathOperator{\supp}{\rm supp\,}

\begin{document}

\title*{Multiplication and linear integral operators on \texorpdfstring{$L_p$}{L_p} spaces representing
polynomial covariance type commutation relations
}
\titlerunning{Multiplication and linear integral operator representing commutation relations }
\author{Domingos Djinja, Sergei Silvestrov, Alex Behakanira Tumwesigye}
\authorrunning{D. Djinja, S. Silvestrov, A. B. Tumwesigye}
\institute{Domingos Djinja \at
Department of Mathematics and Informatics, Faculty of Sciences, Eduardo Mondlane University, Box 257, Maputo, Mozambique. \at
Division of Mathematics and Physics, School of Education, Culture and Communication, M\"alardalens University, Box 883, 72123 V\"aster{\aa}s, Sweden. \\
\email{domingos.djindja@uem.ac.mz, domingos.celso.djinja@mdu.se}
\and
Sergei Silvestrov \at
Division of Mathematics and Physics, School of Education, Culture and Communication, M\"alardalens University,
Box 883, 72123 V\"aster{\aa}s, Sweden. \\
\email{sergei.silvestrov@mdu.se}
\and
Alex Behakanira Tumwesigye \at
Department of Mathematics, College of Natural Sciences, Makerere University, Box 7062, Kampala, Uganda. \\
\email{alexbt@cns.mak.ac.ug}}
%
%


\maketitle
\label{Chap:DjinjaSilvestrovTumwesigye1MultIntoprepscomrels}

\abstract*{Representations of polynomial covariance type commutation relations by pairs of linear integral operators and multiplication operators on Banach spaces $L_p$ are constructed.
}

\abstract{Representations of polynomial covariant type commutation relations by pairs of linear integral operators and multiplication operators on Banach spaces $L_p$ are constructed.   }

\section{Introduction}
\index{Commutation relations!covariance}
Commutation relations of the form
\begin{equation} \label{CovrelABeqBFA}
  AB=B F(A)
\end{equation}
where $A, B$ are elements of an associative algebra and  $F$ is a function of the elements of the algebra, are important in many areas of Mathematics and applications. Such commutation relations are usually called covariance relations, crossed product relations or semi-direct product relations. The pairs $(A,B)$ of elements of an algebra that satisfy \eqref{CovrelABeqBFA} are called representations of this relation
in that algebra. Representations of covariance commutation relations \eqref{CovrelABeqBFA} by linear operators are important for study of actions and induced representations of groups and semigroups, crossed product operator algebras, dynamical systems, harmonic analysis, wavelets and fractals analysis and applications in physics and engineering
\cite{BratJorgIFSAMSmemo99,BratJorgbook,JorgWavSignFracbook,JorgOpRepTh88,JorMoore84,MACbook1,MACbook2,MACbook3,OstSambook,Pedbook79,SamoilenkbookSTDST}.

The structure of representations for the relations of the form \eqref{CovrelABeqBFA}
by bounded and unbounded self-adjoint operators, normal operators, unitary operators, partial isometries
and other linear operators with special involution conditions on a Hilbert space, have been considered in
\cite{BratEvansJorg2000,CarlsenSilvExpoMath07,CarlsenSilvAAM09,CarlsenSilvProcEAS10,DutkayJorg3,DJS12JFASilv,DLS09,DutSilvProcAMS,
DutSilvSV,JSvT12a,JSvT12b,Mansour16,JMusondaPhdth18,JMusonda19,Musonda20,Nazaikinskii96,OstSambook,PerssonSilvestrov031,
PerssonSilvestrov032,PersSilv:CommutRelDinSyst,RST16,RSST16,SamoilenkbookSTDST,SaV8894,SilPhD95,STomdynsystype1,
SilWallin96,SvSJ07a,SvSJ07b,SvSJ07c,SvT09,Tomiyama87,Tomiama:SeoulLN1992,Tomiama:SeoulLN2part2000,TumwesigyePhDThesis2018,TumwRiSilv:ComCrPrAlgPieccnstfnctreallineSPAS19v2,VaislebSa90}
using reordering formulas for functions of the algebra elements and operators satisfying covariance commutation relation,
functional calculus and spectral representation of operators \cite{AkkLnearOperatorsDST} and interplay with dynamical systems generated
by iteration of the maps involved in the commutation relations.

In this paper we construct representations of \eqref{CovrelABeqBFA} by pairs of linear integral and multiplication operators on  Banach space $L_p$.
Such representations can also be viewed as solutions for operator equations $AX=XF(A)$, when $A$ is specified or
$XB=BF(X)$ when $B$ is specified. In contrast to \cite{OstSambook,SamoilenkbookSTDST,SaV8894,VaislebSa90} devoted to involutive representations of covariance type relations by operators on Hilbert spaces using spectral theory of operators on Hilbert spaces,
we aim at direct construction of various classes of representations of covariance type relations in specific important classes of operators on Banach spaces more general than Hilbert spaces without imposing any involution conditions and not using classical spectral theory of operators.
This paper is organized in three sections. After the introduction,  we present in Section \ref{SecPrelNot} preliminaries, notations and basic definitions.
In Section \ref{SecMainResults} we present the main results about construction of specific representations on Banach function spaces $L_p$.

\section{Preliminaries and notations}\label{SecPrelNot}
\index{Banach space}
In this section we present some preliminaries, basic definitions and notations. For more details please read 
\cite{AdamsGDST, DudleyDST,FollandRADST,KantarovitchAkilovFunctionalAnalysis,KolmogorovFominElemTheorFuncFuncAnal,RudinRCADST,RynneLFADST}. \par
Let $S\subseteq \mathbb{R}$, ($\mathbb{R}$ is the set of real numbers), be a Lebesgue measurable set and let $(S,\Sigma, \tilde{m})$ be a $\sigma$-finite measure space, that is, $S$ is a nonempty set, $\Sigma$ is a $\sigma-$algebra with subsets of $S$,
where $S$ can be cover with at most  countable many disjoint sets $E_1,E_2,E_3,\ldots$ such that $ E_i\in \Sigma, \,
\tilde{m}(E_i)<\infty$, $i=1,2,\ldots$  and $\tilde{m}$ is the Lebesgue measure.
For $1\leqslant p<\infty,$ we denote by $L_p(S)$, the set of all classes of equivalent measurable functions $f:S\to \mathbb{R}$ such that
$
    \int\limits_{S} |f(t)|^p dt < \infty.
$
This is a Banach space (Hilbert space when $p=2$) with norm
$ \| f\|_p= \left( \int\limits_{S} |f(t)|^p dt \right)^{\frac{1}{p}}.
$
We denote by $L_\infty(S)$ the set of all classes of equivalent measurable functions $f:S\to \mathbb{R}$ such that there is a constant $\lambda >0$,
$|f(t)|\leq \lambda $
almost everywhere. This is a Banach space with norm
$\displaystyle
  \|f\|_{\infty}=\mathop{\esssup}_{t\in S} |f(t)|.
$

\section{Operator representations of covariance commutation relations}\label{SecMainResults}
Before we proceed with constructions of more complicated operator representations of commutation relations \eqref{CovrelABeqBFA} on more complicated Banach spaces, we wish to
mention the following two observations that, while being elementary, nevertheless explicitly indicate differences in how the different operator representations of commutation relations \eqref{CovrelABeqBFA} interact with the function $F$.
\begin{prop:DjinjaTumwSilv1}
Let $A:E\to E$ and $B:E\to E$, $B\not=0$,  be linear operators on a linear space $E$,
such that any composition among them is well defined and consider $F:\mathbb{R} \to\mathbb{R}$ a polynomial. If $A=\alpha I$, then $AB=BF(A)$ if and only if $F(\alpha)=\alpha$.
\end{prop:DjinjaTumwSilv1}

\begin{proof}
  If $A=\alpha I$, then
  \begin{gather*}
    AB=\alpha IB=\alpha B, \\
    BF(A)=BF(\alpha I)=BF(\alpha) I=F(\alpha)B.
  \end{gather*}
  We have then $AB=BF(A)$, $B\neq 0$ if and only if $F(\alpha)=\alpha$.
\qed \end{proof}

\begin{prop:DjinjaTumwSilv1}
Let $A:E\to E$ and $B:E\to E$ be linear operators such that any composition among them is well defined and consider a polynomial $F:\mathbb{R} \to\mathbb{R}$. If $B=\alpha I$, where $\alpha\not=0$,  then $AB=BF(A)$ if and only if $F$ satisfies $F(A)=A$.
\end{prop:DjinjaTumwSilv1}

\begin{proof}
  If $B=\alpha I$ then
  \begin{gather*}
    AB= A(\alpha I)=\alpha A, \\
    BF(A)=\alpha IF(A)=\alpha F(A) .
  \end{gather*}
  We have then $AB=BF(A)$ if and only if $F(A)=A$.
\qed \end{proof}

\subsection{Representations of covariance commutation relations by integral and multiplication operators on \texorpdfstring{$L_p$}{L_p} spaces}
We consider first a useful lemma for integral operators.
\begin{lem:DjinjaTumwSilv1}\label{lemEqIntForAllLpFunct}
Let $f:[\alpha_1,\beta_1]\to \mathbb{R}$,  $g: [\alpha_2,\beta_2]\to \mathbb{R}$ be two measurable functions such that
for all $x\in L_p(\mathbb{R})$, $1\leq p \leq \infty $,
\begin{equation*}
  \int\limits_{\alpha_1}^{\beta_1} f(t)x(t)dt,\quad \int\limits_{\alpha_2}^{\beta_2} g(t)x(t)dt,
\end{equation*}
exist and are finite, where  $\alpha_1,\beta_1,\alpha_2,\beta_2\in \mathbb{R}$, $\alpha_1<\beta_1$ and $\alpha_2<\beta_2$.
Set
$
  G=[\alpha_1,\beta_1]\cap [\alpha_2,\beta_2].
$
 Then the following statements are equivalent.
\begin{enumerate}
  \item For all $x\in L_p(\mathbb{R})$, where  $1\le p\le \infty$,  the following holds
      \begin{equation*} \label{EqInteDiffIntervalsLp}
  \int\limits_{\alpha_1}^{\beta_1} f(t)x(t)dt=\int\limits_{\alpha_2}^{\beta_2} g(t)x(t)dt.
\end{equation*}
\item The following conditions hold:
    \begin{enumerate}
      \item[a)] for almost every $t\in G$, $f(t)=g(t)$;
      \item[b)] for almost every $t \in [\alpha_1,\beta_1]\setminus G,\ f(t)=0;$
      \item[c)]  for almost every $t \in [\alpha_2,\beta_2]\setminus G,\ g(t)=0.$
    \end{enumerate}
\end{enumerate}
\end{lem:DjinjaTumwSilv1}

\begin{proof}
    $2.\, \Rightarrow \, 1.$ follows from direct computation. \par
     Suppose that 1. is true.
     Take
    $x(t)=I_{G_1}(t)$ the indicator function of the set $G_1=[\alpha_1,\beta_1]\cup[\alpha_2,\beta_2]$. For this
    function we have,
    \begin{gather*}
      \int\limits_{\alpha_1}^{\beta_1} f(t)x(t)dt=\int\limits_{\alpha_2}^{\beta_2} g(t)x(t)dt=
      \int\limits_{\alpha_1}^{\beta_1} f (t)dt=\int\limits_{\alpha_2}^{\beta_2} g (t)dt =\eta,
    \end{gather*}
    $\eta$ is a constant. Now by taking $x(t)=I_{[\alpha_1,\beta_1]\setminus G}(t)$ we get
         \begin{gather*}
      \int\limits_{\alpha_1}^{\beta_1} f(t)x(t)dt=\int\limits_{\alpha_2}^{\beta_2} g(t)x(t)dt=
      \int\limits_{[\alpha_1,\beta_1]\setminus G} f (t)dt=\int\limits_{\alpha_2}^{\beta_2} g (t)\cdot 0dt =0.
    \end{gather*}
    Then
    $
      \int\limits_{[\alpha_1,\beta_1]\setminus G} f (t)dt=0.
    $
    If instead $x(t)=I_{[\alpha_2,\beta_2]\setminus G}(t)$, then
    $ \int\limits_{[\alpha_2,\beta_2]\setminus G} g (t)dt=0.
    $
    We claim that $f(t)=0$ for almost every $t\in [\alpha_1,\beta_1]\setminus G$ and
    $g(t)=0$ for almost every $t\in [\alpha_2,\beta_2]\setminus G$. We take a partition $\bigcup S_i$ of the set  $[\alpha_1,\beta_1]\setminus G$
    such that $S_i\cap S_j\neq \emptyset$ for $i\neq j$ and each set $S_i$ has positive measure. For each $x(t)=I_{S_i}(t)$, we have
    \begin{gather*}
    \int\limits_{\alpha_1}^{\beta_1} f(t)x(t)dt=\int\limits_{\alpha_2}^{\beta_2} g(t)x(t)dt=
    \int\limits_{S_i} f (t)dt=\int\limits_{\alpha_2}^{\beta_2} g (t)\cdot 0dt =0.
    \end{gather*}
    Thus,
    $\int\limits_{S_i} f (t)dt=0.
    $
    Since we can choose arbitrary partition with positive measure on each of its elements we have
    \begin{gather*}
      f(t)=0 \quad \mbox{ for almost every } t\in [\alpha_1,\beta_1]\setminus G.
    \end{gather*}
    Analogously,
    $g(t)=0$ for almost every $t\in [\alpha_2,\beta_2]\setminus G.$
    Then,
    $$
      \eta = \int\limits_{\alpha_1}^{\beta_1} f (t)dt=\int\limits_{\alpha_2}^{\beta_2} g (t)dt =\int\limits_G f(t)dt=\int\limits_G g (t)dt.
    $$
    Then, for all function $x\in L_p(\mathbb{R})$ we have
    \begin{gather*}
      \int\limits_G f (t)x(t)dt=\int\limits_G g (t)x(t)dt \Longleftrightarrow   \int\limits_G [f(t)-g(t)]x(t)dt=0.
    \end{gather*}
    By taking $x(t)=\left\{\begin{array}{cc} 1, & \mbox{if } f(t)-g(t)>0, \\ -1,  & \mbox{ if } f(t)-g(t)<0,\,   \end{array}\right.$
    for almost every $t\in G$ and $x(t)=0$ for almost every $t\in \mathbb{R}\setminus G$, we get
    $\int\limits_{G} |f(t)-g(t)|dt =0$. This implies that $f(t)=g(t)$ for almost every $t\in G$.
\qed \end{proof}

\begin{rem:DjinjaTumwSilv1}{\rm
   When operators are given in abstract form,  we use the notation $A:L_p(\mathbb{R})\to L_p(\mathbb{R})$ meaning that operator $A$  is well defined
   from $L_p(\mathbb{R})$ to $L_p(\mathbb{R})$ without discussing sufficient conditions for it to be satisfied. For instance, for the following integral operator
   \begin{equation*}
    (Ax)(t)= \int\limits_{\mathbb{R}} k(t,s)x(s)ds
   \end{equation*}
    there are sufficient conditions on the kernel $k(\cdot,\cdot)$ such that operator $A$ is well defined as a linear operator
    from $L_p(\mathbb{R}) $ to $L_p(\mathbb{R})$, and bounded  \cite{ConwayFunctionalAnalysisDST,FollandRADST}. For instance,  \cite[Theorem 6.18]{FollandRADST} states the following:
   if $1\leq p \leq \infty$ and  $k:\mathbb{R}\times [\alpha,\beta]\to \mathbb{R}$  is a measurable function, $\alpha,\beta \in \mathbb{R}$, $\alpha<\beta$, and there is a constant $\lambda>0$ such that
\begin{equation*}
   \mathop{\esssup}_{s\in [\alpha,\beta]}\int\limits_{\mathbb{R}}|k(t,s)| dt\leq \lambda, \quad \mathop{\esssup}_{t\in \mathbb{R}}\int\limits_{\alpha}^\beta|k(t,s)| ds \leq \lambda,
\end{equation*}
 then $A$ is well defined from $L_p(\mathbb{R})$ to $L_p(\mathbb{R})$, $1\leq p \leq \infty$ and bounded.
}
 \end{rem:DjinjaTumwSilv1}

\subsubsection{Representations when \texorpdfstring{$A$}{A} is integral operator and \texorpdfstring{$B$}{B} is multiplication operator}

\begin{prop:DjinjaTumwSilv1}\label{propIntOpLpAintBidentity}       \index{Operator!integral}
Let $A:L_p(\mathbb{R})\to L_p(\mathbb{R})$,  $B:L_p(\mathbb{R})\to L_p(\mathbb{R})$, $1\leq p \leq \infty$, be defined as follows, for almost all $t\in \mathbb{R}$,
\begin{gather*}
  (Ax)(t)= \int\limits_{\alpha}^{\beta}k(t,s) x(s)ds,\quad (Bx)(t)= b(t) x(t), \quad \alpha,\beta\in \mathbb{R}, \ \alpha<\beta,
\end{gather*}
where $k:\mathbb{R}\times [\alpha,\beta]\to \mathbb{R}$ is a measurable function,
and $b:\mathbb{R}\to \mathbb{R}$ is a measurable function.
Consider a polynomial $F(z)=\sum\limits_{j=0}^{n}\delta_j z^j$, where $\delta_0,\ldots,\delta_n\in \mathbb{R}$.
We set
\begin{gather} \nonumber
  k_0(t,s)=k(t,s), \quad 
 \quad k_m(t,s)=\int\limits_{\alpha}^{\beta} k(t,\tau)k_{m-1}(\tau,s)d\tau,\quad m\in \{1,\ldots,n\}
  \\
F_0(k(t,s))=0,\quad F_n(k(t,s))=\sum_{j=1}^{n} \delta_j k_{j-1}(t,s),\quad n\geq 1.
\label{FnkernelForFAiteredKernerlns}
\end{gather}
Then $ AB=BF(A)$
if and only if
\begin{equation}\label{CondABeqBFADeltaZeroIntOpIdentity}
\forall \ x\in L_p(\mathbb{R}): \quad  b(t)\delta_0 x(t)+
   b(t)\int\limits_{\alpha}^{\beta} F_n(k(t,s))x(s)ds = \int\limits_{\alpha}^{\beta}  k(t,s)b(s) x(s)ds.
\end{equation}
If $\delta_0=0$, that is, $F(z)=\sum\limits_{j=1}^{n}\delta_j z^j$, then
the condition \eqref{CondABeqBFADeltaZeroIntOpIdentity} reduces to the following: for almost every
 $(t,s)$ in $\mathbb{R}\times [\alpha,\beta]$,
\begin{gather}\label{CondABeqBFADeltaZeroEqualZeroIntOpIndet}
  b(t)F_n(k(t,s))=k(t,s)b(s).
\end{gather}
\end{prop:DjinjaTumwSilv1}

\begin{proof}
By applying Fubini Theorem from \cite{AdamsGDST} and iterative kernels from \cite{KrasnolskZabreykoPuSoIntegralopssummablefuncs}, We have
  \begin{eqnarray*}
  (A^2x)(t)&=&\int\limits_{\alpha}^{\beta} k(t,s)(Ax)(s)ds=
  \int\limits_{\alpha}^{\beta} k(t,s)\left(\int\limits_{\alpha}^{\beta} k(s,\tau)x(\tau)d\tau\right)ds\\
  &=&\int\limits_{\alpha}^{\beta} \left(\int\limits_{\alpha}^{\beta} k(t,s)k(s,\tau)ds\right)x(\tau) d\tau =
  \int\limits_{\alpha}^{\beta} k_1(t,\tau)x(\tau)d\tau,
  \end{eqnarray*}
  where
  $
    k_1(t,s)=\int\limits_{\alpha}^{\beta} k(t,\tau)k(\tau,s)d\tau.
  $
 In the same way,
  \begin{gather*}
    (A^3x)(t)=\int\limits_{\alpha}^{\beta} k(t,s)(A^2x)(s)ds
  =\int\limits_{\alpha}^{\beta} k(t,s)\left(\int\limits_{\alpha}^{\beta}  k_1(s,\tau)x(\tau)d\tau\right)ds \\
  =\int\limits_{\alpha}^{\beta} k_2(t,s)x(s)ds,
  \end{gather*}
  where $k_2(t,s)=\int\limits_{\alpha}^{\beta} k(t,\tau)k_1(\tau,s)d\tau.$
For every $m \ge 1$,
$$
  (A^m x)(t)=\int\limits_{\alpha}^{\beta} k_{m-1}(t,s)x(s)ds,
$$
where $k_m(t,s)=\int\limits_{\alpha}^{\beta} k(t,\tau)k_{m-1}(\tau,s)d\tau,\ m=1,\ldots,n,\ k_0(t,s)=k(t,s).$
Thus,
 \begin{eqnarray*}
   (F(A)x)(t)&=&\delta_0 x(t)+ \sum\limits_{j=1}^{n} \delta_j (A^j x)(t)=\delta_0 x(t)+\sum\limits_{j=1}^{n} \delta_{j} \int\limits_{\alpha}^{\beta} k_{j-1}(t,s)x(s)ds \\
   &=&\delta_0 x(t)+\int\limits_{\alpha}^{\beta} F_n(k(t,s))x(s)ds,
   \end{eqnarray*}
  where
  $
      F_n(k(t,s))=\sum\limits_{j=1}^{n} \delta_j k_{j-1}(t,s),
  $
for $n\geq 1,$ and $F_0(k(t,s))=0.$
  So, we can compute $BF(A)x$ and $(AB)x$   as follows:
   \begin{align*}
     (BF(A)x)(t) & =b(t)(F(A)x)(t)=b(t)\delta_0 x(t)+
   b(t)\int\limits_{\alpha}^{\beta} F_n(k(t,s))x(s)ds,   \\
   (ABx)(t) & =A(Bx)(t)= \int\limits_{\alpha}^{\beta}  k(t,s)b(s) x(s)ds.
  \end{align*}
It follows that $ABx=BF(A)x$ if and only if condition \eqref{CondABeqBFADeltaZeroIntOpIdentity} holds.

If  $\delta_0=0$ then  condition \eqref{CondABeqBFADeltaZeroIntOpIdentity} reduces to the following:
\begin{gather*}
 \forall \ x\in L_p(\mathbb{R}): \quad   \int\limits_{\alpha}^{\beta} b(t)F_n(k(t,s))x(s)ds = \int\limits_{\alpha}^{\beta}  k(t,s)b(s) x(s)ds.
\end{gather*}
Let $f(t,s)= b(t)F_n(k(t,s))-k(t,s)b(s)$. By applying Lemma \ref{lemEqIntForAllLpFunct}
we have for
almost every $t\in \mathbb{R}$ that $f(t,\cdot)=0$ almost everywhere.  
Since the set $N=\{(t,s)\in\mathbb{R}\times[\alpha,\beta]:\ f(t,s)\not=0\}\subset \mathbb{R}^2$ is measurable
and  almost all sections $N_t=\{s\in [\alpha,\beta]:\ (t,s)\in N \}$ of the plane has Lebesgue measure zero, by
the reciprocal Fubini Theorem \cite{OxtobyDST}, the set $N$ has Lebesgue measure zero on the plane $\mathbb{R}^2$.
\qed \end{proof}

\begin{cor:DjinjaTumwSilv1}\label{CorIntOpMultOpPolinomials}  For $M_1,M_2\in\mathbb{R}$, $M_1<M_2$ and  $1\leq p\leq\infty$, let $A:L_p([M_1,M_2])\to L_p([M_1,M_2])$ and $B:L_p([M_1,M_2])\to L_p([M_1,M_2])$  be nonzero operators defined, for almost all $t$, by
\begin{gather*}
  (Ax)(t)= \int\limits_{\alpha}^{\beta}k(t,s) x(s)ds,\quad (Bx)(t)= b(t) x(t), \quad \alpha,\beta\in \mathbb{R},\ \alpha<\beta,
\end{gather*}
where $[M_1,M_2]\supseteq [\alpha,\beta]$, and $k(\cdot,\cdot):[M_1,M_2]\times [\alpha,\beta]\to \mathbb{R}$, $b:[M_1,M_2]\to\mathbb{R}$ are given by
\begin{gather*}
  k(t,s)=a_0+a_1t+c_1s,\qquad b(t)=\sum\limits_{j=0}^n b_j t^j,
\end{gather*}
where $n$ is non-negative integer, $a_0,\, a_1,\, c_1,\, b_j $ are real numbers for $j=0,\ldots,n$. Consider a polynomial  defined by $F(z)=\delta_0+\delta_1 z +\delta_2 z^2$, where $\delta_0,\,\delta_1,\, \delta_2\in\mathbb{R}$.

Then, $AB=BF(A)$ if and only if
\begin{gather*}
\forall \ x\in L_p([M_1,M_2]): \quad  b(t)\delta_0 x(t)+
   b(t)\int\limits_{\alpha}^{\beta} F_n(k(t,s))x(s)ds = \int\limits_{\alpha}^{\beta}  k(t,s)b(s) x(s)ds,
\end{gather*}
where $F_n(k(t,s))$ is given by \eqref{FnkernelForFAiteredKernerlns}.

If $\delta_0=0$, that is, $F(z)=\delta_1 z+\delta_2 z^2$ then
the last condition reduces to the condition that for almost every
 $(t,s)$ in $[M_1,M_2]\times [\alpha,\beta]$
\begin{equation}\label{CondCommuationRelationCorPolyn}
  b(t)F_2(k(t,s))=k(t,s)b(s).
\end{equation}
Condition \eqref{CondCommuationRelationCorPolyn} is equivalent to that $b(\cdot)\equiv b_0 \neq 0$ is a nonzero constant
$(b_j=0,$ $j=1,\ldots,n)$ and one of the following cases holds:
\begin{enumerate}
  \item\label{corIntOpMultOpPolynomialb0:item1} if $\delta_2=0$, $\delta_1=1$, then $a_0,a_1,c_1\in\mathbb{R}$ can be arbitrary;

  \item\label{corIntOpMultOpPolynomialb0:item2} if $\delta_2\not=0$, $\delta_1=1$, $a_1\not=0$, $c_1=0$, then
  \begin{equation*}
    a_0=-\frac{\beta+\alpha}{2}a_1;
  \end{equation*}

    \item\label{corIntOpMultOpPolynomialb0:item3} if $\delta_2\not=0$, $\delta_1=1$, $a_1=0$, $c_1\not=0$, then
  \begin{equation*}
    a_0=-\frac{\beta+\alpha}{2}c_1;
  \end{equation*}

  \item\label{corIntOpMultOpPolynomialb0:item4} if $\delta_2\not=0$, $\delta_1\not=1$, $a_1\not=0$, $c_1=0$,  then
  \begin{equation*}
          a_0=\frac{2-2\delta_1-\delta_2(\beta^2-\alpha^2)a_1}{2\delta_2(\beta-\alpha)};
  \end{equation*}

  \item\label{corIntOpMultOpPolynomialb0:item5} if $\delta_2\not=0$, $\delta_1\not=1$, $c_1\not=0$, $a_1=0$, then
           \begin{equation*}
          a_0=\frac{2-2\delta_1-\delta_2(\beta^2-\alpha^2)c_1}{2\delta_2(\beta-\alpha)};
        \end{equation*}

   \item\label{corIntOpMultOpPolynomialb0:item6} if $\delta_2\not=0$, $\delta_1\not=1$, $a_1=0$ and $c_1=0$,  then
  \begin{equation*}
    a_0=\frac{1-\delta_1}{\delta_2(\beta-\alpha)}.
  \end{equation*}
\end{enumerate}
\end{cor:DjinjaTumwSilv1}

\begin{proof}
  Operator  $A$ is well defined on $L_p[M_1,M_2]$, $1\leq p \leq \infty $. In fact, it follows by \cite[Theorem 3.4.10]{HutsonPymDST}.
  Moreover, kernel $k(\cdot,\cdot)$  is continuous on a closed and bounded set
  $[-M,M]\times [\alpha,\beta]$ and $b(\cdot)$ is continuous in $[M_1,M_2]$, so these functions are measurable.
  By applying Proposition \ref{propIntOpLpAintBidentity} we just need to check when the  condition \eqref{CondABeqBFADeltaZeroEqualZeroIntOpIndet} is satisfied for $k(\cdot,\cdot)$ and $b(\cdot)$. 
  We compute
\begin{align}\nonumber
k_1(t,s)=&\int\limits_{\alpha}^{\beta}k(t,\tau)k(\tau,s)d\tau =\int\limits_{\alpha}^{\beta}(a_0+a_1t+c_1\tau)(a_0+a_1\tau+c_1s)d\tau\\ \nonumber
 =& \int\limits_{\alpha}^{\beta} [(a_0^2+a_0a_1t+a_0c_1s+a_1c_1ts)\\  \nonumber &+(a_0a_1+a_0c_1+a_1^2t+c_1^2s)\tau+a_1c_1\tau^2]d\tau
 \\ \nonumber
 =& (\beta-\alpha)(a_0^2+a_0a_1t+a_0c_1s+a_1c_1ts)\\ \nonumber
 &+\frac{\beta^2-\alpha^2}{2}\cdot (a_0a_1+a_0c_1+a_1^2t+c_1^2s)\\
 &+\frac{\beta^3-\alpha^3}{3}a_1c_1
=\nu_0 + \nu_1 t +\nu_2 s +\nu_3 ts, \label{KOneRemCorAIntOpBIndentOP}
\end{align}
where
\begin{gather*}
\begin{array}{cclccl}
  \nu_0&=&a_0^2(\beta-\alpha)+\frac{\beta^2-\alpha^2}{2}a_0(a_1+c_1)+a_1c_1 \frac{\beta^3-\alpha^3}{3},\hspace{0mm} & \nu_2&=&a_0 c_1(\beta-\alpha)+c_1^2 \frac{\beta^2-\alpha^2}{2},  \\
  \nu_1&=&a_1^2\frac{\beta^2-\alpha^2}{2}+a_1a_0 (\beta-\alpha), \hspace{-4mm}& \nu_3&=&a_1c_1(\beta-\alpha).
\end{array}
\end{gather*}
Then, we have
\begin{align*}
& b(t)F_2(k(t,s))=b(t)(\delta_1 k(t,s)+\delta_2 k_1(t,s) )=(a_0 \delta_1 +\delta_2 \nu_0)\sum_{j=0}^n b_j t^j  \\
               &+(a_1\delta_1+\delta_2 \nu_1) \sum_{j=0}^n b_j t^{j+1}+
              ( c_1\delta_1+\delta_2 \nu_2) \sum_{j=0}^n b_j t^js
              + \nu_3 \delta_2 \sum_{j=0}^n b_j t^{j+1}s \\
         &\resizebox{0.99\hsize}{!}{$\displaystyle= (\delta_1a_0+\delta_2\nu_0)b_0+(c_1\delta_1+\nu_2\delta_2)b_0s+
         \sum_{j=1}^{n}((\delta_1a_0+\delta_2\nu_0)b_j+(\delta_1a_1+\delta_2\nu_1)b_{j-1})t^j
         $}
         \\
         &+
         \sum_{j=1}^{n}((c_1\delta_1+\nu_2\delta_2)b_j+\nu_3\delta_2b_{j-1})t^js+(\delta_1 a_1+\delta_2\nu_1)b_nt^{n+1}+\nu_3\delta_2 b_n t^{n+1}s\\
& k(t,s)b(s)=a_0\sum_{j=0}^n b_j s^j+a_1 \sum_{j=0}^n b_j s^j t +c_1 \sum_{j=0}^n b_j s^{j+1}=a_0b_0+a_1b_0t
\\
&
+\sum_{j=1}^{n} (a_0b_j+c_1b_{j-1})s^j 
+\sum_{j=1}^{n}a_1b_js^j t + c_1b_ns^{n+1}.
\end{align*}
Thus we have $\displaystyle k(t,s)b(s)=b(t)F_2(k(t,s))$ for all $(t,s)\in [M_1,M_2]\times[\alpha,\beta]$ if and only if
 \begin{eqnarray} \nonumber
 a_0b_0&=& (a_0\delta_1  + \delta_2 \nu_0 )b_0 \\ \nonumber
 a_1b_0&=&(a_0\delta_1+\delta_2 \nu_0 )b_1+(a_1\delta_1  +\delta_2 \nu_1) b_0\\ \label{EqCorProofb1=0IntOpMutOpPolynomials}
 a_0b_1+c_1b_0&=& (c_1\delta_1+\delta_2 \nu_2) b_0\\ \label{EqCorProofb1=0a_1=0IntOpMutOpPolynomials}
 a_1b_1 &=& (c_1\delta_1  + \delta_2 \nu_2) b_1+ \delta_2 \nu_3 b_0 \\
 \label{EqSevenRemCorIntOPMultOpPolinomials}
 0&=& a_0b_j +c_1 b_{j-1},\quad 2\le j\le n \\ \nonumber
 0&=&  (a_0\delta_1+\delta_2 \nu_0) b_j+(a_1\delta_1 +\delta_2 \nu_1) b_{j-1},\quad 2\le j\le n\\
\label{EqElevenRemCorIntOpMultOpPolinomials}
 a_1b_j&=& 0,\quad 2\le j \le n\\
 \nonumber
 0&=& c_1\delta_1b_j+\delta_2 \nu_3 b_{j-1}+\delta_2 \nu_2 b_j\quad 2\le j\le n\\
 \nonumber
0&=& a_1\delta_1 b_n+\delta_2 \nu_1 b_n,\quad \mbox{ if }  n\geq 1 \\  \label{Eqc1bnLastCoefPolynomial}
 c_1b_{n}&=&0,\quad \mbox{ if } n\ge 1 \\  \nonumber
 0&=& \delta_2 \nu_3 b_n,\quad \mbox{ if } n\ge 1.
 \end{eqnarray}

Suppose that $n\ge 1$. We proceed by induction to prove that $b_j=0$, for all $j=1,\ldots,n$. For $i=0$, we suppose that $b_n=b_{n-i}\not=0$. Then
from the equation \eqref{EqElevenRemCorIntOpMultOpPolinomials} we have  $a_1b_n=0$ and thus $a_1=0$. From the equation \eqref{Eqc1bnLastCoefPolynomial} we have $c_1b_n=0$ and thus $c_1=0$. From the equation
\eqref{EqSevenRemCorIntOPMultOpPolinomials} we have $0= a_0b_n +c_1 b_{n-1}=a_0b_n$ and thus $a_0=0$. This implies that $k(t,s)\equiv 0$, that is, $A=0$. So for $i=0$, $b_n=b_{n-i}\not= 0$ implies
$A= 0$. Hence,  $b_n=0$.
Let $1<m\leq n-2$ and suppose that $b_{n-i}=0$ for all $i=1,\ldots,m-1$. Let us show that then $b_{n-m}=0$.  If $b_{n-m}\not=0$, then from equation \eqref{EqElevenRemCorIntOpMultOpPolinomials} we have $a_1b_{n-m}=0$ which implies $a_1=0$. From equation
\eqref{EqSevenRemCorIntOPMultOpPolinomials} and for $j=n-m+1$ by induction assumption  $a_0b_{n-m+1}+c_1b_{n-m}=c_1b_{n-m}=0$ which implies
$c_1=0$. Therefore  from \eqref{EqSevenRemCorIntOPMultOpPolinomials} and for $j=n-m$ we have $a_0b_{n-m}=0$
which implies $a_0=0$. Then $k(t,s)\equiv 0$, that is $A= 0$. So we must have $b_{n-m}=0$.
If $m=n-1$, then let us show that $b_{n-m}=b_1=0$. If $b_{n-m}\not=0$
then \eqref{EqSevenRemCorIntOPMultOpPolinomials} gives $c_1b_{n-m}=c_1b_1=0$ when $j=n-m+1=2$.
Then $c_1=0$ and by \eqref{EqCorProofb1=0a_1=0IntOpMutOpPolynomials}, since $\nu_2=\nu_3=0$  we get $a_1b_1=0$
which yields $a_1=0$. Therefore, \eqref{EqCorProofb1=0IntOpMutOpPolynomials} gives $a_0b_1=0$
which yields $a_0=0$. Thus $A=0$. Since $A\not=0$, $b_1=0$ is proved.
 Thus $b(\cdot)=b_0$ is proved.

 Since $B\neq 0$ and $B=b_0 I$ (multiple of identity operator),  $b_0\neq 0$ and the commutation relation is equivalent to $F(A)=A$.  By \eqref{CondABeqBFADeltaZeroEqualZeroIntOpIndet}  we have
$F_2(k(t,s))=k(t,s)$ which can be written as follows
\begin{equation}\label{EqF(A)=AintermsOfKernelsCorIntOpMultOpPolynomial}
  \delta_1 k(t,s)+\delta_2 k_1(t,s)=k(t,s),
\end{equation}
where $k(t,s)=a_0+a_1t+c_1s$ and $k_1(t,s)=\nu_0+\nu_1 t+\nu_2 s+\nu_3 ts$,
\begin{gather*}
\begin{array}{cclccl}
  \nu_0&=&a_0^2(\beta-\alpha)+\frac{\beta^2-\alpha^2}{2}a_0(a_1+c_1)+a_1c_1 \frac{\beta^3-\alpha^3}{3},\hspace{0mm} & \nu_2&=&a_0 c_1(\beta-\alpha)+c_1^2 \frac{\beta^2-\alpha^2}{2},  \\
  \nu_1&=&a_1^2\frac{\beta^2-\alpha^2}{2}+a_1a_0 (\beta-\alpha), \hspace{-4mm}& \nu_3&=&a_1c_1(\beta-\alpha).
\end{array}
\end{gather*}

If $\delta_2=0$, then  \eqref{EqF(A)=AintermsOfKernelsCorIntOpMultOpPolynomial} becomes $(\delta_1-1)k(\cdot,\cdot)=0$ and $A\not=0$ yields $\delta_1=1$. Thus, if $\delta_2=0$ and $\delta_1=1$, then \eqref{EqF(A)=AintermsOfKernelsCorIntOpMultOpPolynomial} is satisfied for any $a_0,a_1,c_1\in\mathbb{R}$.

 If $\delta_2\not=0$ and $\delta_1=1$ then \eqref{EqF(A)=AintermsOfKernelsCorIntOpMultOpPolynomial}
 becomes $k_1(\cdot,\cdot)=0$, that is, $\nu_0=\nu_1=\nu_2=\nu_3=0$, where
 \begin{gather*}
\begin{array}{cclccl}
  \nu_0&=&a_0^2(\beta-\alpha)+\frac{\beta^2-\alpha^2}{2}a_0(a_1+c_1)+a_1c_1 \frac{\beta^3-\alpha^3}{3},\hspace{0mm} & \nu_2&=&a_0 c_1(\beta-\alpha)+c_1^2 \frac{\beta^2-\alpha^2}{2},  \\
  \nu_1&=&a_1^2\frac{\beta^2-\alpha^2}{2}+a_1a_0 (\beta-\alpha), \hspace{-4mm}& \nu_3&=&a_1c_1(\beta-\alpha).
\end{array}
\end{gather*}
 Since $\alpha<\beta$, $a_1c_1(\beta-\alpha)=0$ is equivalent to either $a_1=0$ or $c_1=0$.  If $a_1\neq0$, $c_1=0$, then
 \begin{gather*}
   \left\{\begin{array}{c}
            \nu_0=0  \\
            \nu_1=0\\
            \nu_2=0 \\
            \nu_3=0
          \end{array}\right. \Leftrightarrow
             \left\{\begin{array}{c}
            (\beta-\alpha)a_0^2+\frac{\beta^2-\alpha^2}{2}a_0a_1=0  \\ \medskip
            (\beta-\alpha)a_1a_0+\frac{\beta^2-\alpha^2}{2}a^2_1=0
          \end{array}\right. \Leftrightarrow a_0+\frac{\beta+\alpha}{2}a_1=0,
 \end{gather*}
 which is equivalent to $a_0=-\frac{\beta+\alpha}{2}a_1$. If $a_1=0$, $c_1\neq0$, then
 \begin{gather*}
   \left\{\begin{array}{c}
            \nu_0=0  \\
            \nu_1=0\\
            \nu_2=0 \\
            \nu_3=0
          \end{array}\right. \Leftrightarrow
             \left\{\begin{array}{c}
            (\beta-\alpha)a_0^2+\frac{\beta^2-\alpha^2}{2}a_0c_1=0  \\ \medskip
            (\beta-\alpha)c_1a_0+\frac{\beta^2-\alpha^2}{2}c^2_1=0
          \end{array}\right. \Leftrightarrow a_0+\frac{\beta+\alpha}{2}c_1=0,
 \end{gather*}
 which is equivalent to $a_0=-\frac{\beta+\alpha}{2}c_1$. If $a_1=0$, $c_1=0$, then
 $\nu_0=\nu_1=\nu_2=\nu_3=0$ is equivalent to $a_0^2(\beta-\alpha)=0$, that is, $a_0=0$. This implies $A=0$. Therefore, $\delta_2\not=0$, $\delta_1=1$, $a_1=c_1=0$ yields $A=0$.

Consider $\delta_2\not=0$ and $\delta_1\not=1$, and note that  \eqref{EqF(A)=AintermsOfKernelsCorIntOpMultOpPolynomial} is equivalent to:
 \begin{gather}\label{SystEqCorIntOpMultOpPolynomialdelta0zerodelta2not0}
   \left\{\begin{array}{cl}
          a_0&= \delta_1a_0+ \delta_2 a_0^2(\beta-\alpha)+\delta_2\frac{\beta^2-\alpha^2}{2}a_0(a_1+c_1)+\delta_2a_1c_1 \frac{\beta^3-\alpha^3}{3}  \\
          a_1&= \delta_1a_1+\delta_2 a_1^2\frac{\beta^2-\alpha^2}{2}+\delta_2 a_1a_0 (\beta-\alpha)\\
           c_1&= \delta_1 c_1+\delta_2 a_0 c_1(\beta-\alpha)+\delta_2 c_1^2 \frac{\beta^2-\alpha^2}{2} \\
           0&= \delta_2 a_1c_1(\beta-\alpha).
          \end{array}\right.
 \end{gather}
Since $\alpha <\beta$ and $\delta_2\not=0$,  equation $\delta_2 a_1c_1(\beta-\alpha)=0$ implies that either $a_1=0$ or $c_1=0$.
If $\delta_2\not=0$, $\delta_1\not=1$, $a_1\not=0$ and $c_1=0$, then  \eqref{SystEqCorIntOpMultOpPolynomialdelta0zerodelta2not0} becomes
 \begin{eqnarray*}
a_0&=& \delta_1a_0+ \delta_2 a_0^2(\beta-\alpha)+\delta_2\frac{\beta^2-\alpha^2}{2}a_0a_1\\
          a_1&=& \delta_1a_1+\delta_2 a_1^2\frac{\beta^2-\alpha^2}{2}+\delta_2 a_1a_0 (\beta-\alpha)
 \end{eqnarray*}
 which is equivalent to $1=\delta_1+\delta_2(\beta-\alpha)a_0+\delta_2 \frac{\beta^2-\alpha^2}{2}a_1$. Then,
 \begin{equation*}
   a_0=\frac{2-2\delta_1-\delta_2(\beta^2-\alpha^2)a_1}{2\delta_2(\beta-\alpha)}.
 \end{equation*}

If $\delta_2\not=0$, $\delta_1\not=1$, $a_1=0$ and $c_1\neq0$, then \eqref{SystEqCorIntOpMultOpPolynomialdelta0zerodelta2not0} becomes
 \begin{eqnarray*}
a_0&=& \delta_1a_0+ \delta_2 a_0^2(\beta-\alpha)+\delta_2\frac{\beta^2-\alpha^2}{2}a_0c_1\\
          c_1&=& \delta_1c_1+\delta_2 c_1^2\frac{\beta^2-\alpha^2}{2}+\delta_2 c_1a_0 (\beta-\alpha)
 \end{eqnarray*}
 which is equivalent to $1=\delta_1+\delta_2(\beta-\alpha)a_0+\delta_2 \frac{\beta^2-\alpha^2}{2}c_1$. Then,
 \begin{equation*}
   a_0=\frac{2-2\delta_1-\delta_2(\beta^2-\alpha^2)c_1}{2\delta_2(\beta-\alpha)}.
 \end{equation*}

If $\delta_2\not=0$, $\delta_1\not=1$, $a_1=0$ and $c_1=0$, then  $A\not=0$ yields $a_0\not=0$ and \eqref{SystEqCorIntOpMultOpPolynomialdelta0zerodelta2not0} becomes
 \begin{eqnarray*}
&&           a_0=\delta_1a_0+ \delta_2 a_0^2(\beta-\alpha)
 \end{eqnarray*}
 which is equivalent to
 $
   a_0=\frac{1-\delta_1}{\delta_2(\beta-\alpha)}.
 $
 \qed
 \end{proof}

\begin{rem:DjinjaTumwSilv1} {\rm
The integral operator given by $(Ax)(t)=\int\limits_{\alpha_1}^{\beta_1} k(t,s)x(s)ds$  for almost all $t$, where $k:[\alpha_1,\beta_1]\times[\alpha_1,\beta_1]\to \mathbb{R}$ is a measurable function that satisfies
\begin{equation*}
  \int\limits_{\alpha_1}^{\beta_1} \left(\int\limits_{\alpha_1}^{\beta_1} |k(t,s)|^qds\right)^\frac{p}{q}dt <\infty,
\end{equation*}
by {\rm \cite[Theorem 3.4.10]{HutsonPymDST}} is a well defined and bounded linear operator
from $L_p{[\alpha_1,\beta_1]}$ to $L_p{[\alpha_1 ,\beta_1]}$, where  $1<p<\infty$, $1<q<\infty$, $\frac{1}{p}+\frac{1}{q}=1$.
}
\end{rem:DjinjaTumwSilv1}

\begin{rem:DjinjaTumwSilv1}{\rm
If in the Corollary \ref{CorIntOpMultOpPolinomials} when $0\not\in [M_1,M_2]$, one  takes $b(t)$ to be a Laurent polynomial with only negative powers of $t$ then there is no  non-zero
kernel $k(t,s)=a_0+a_1t+c_1s$ (there is no $A\not=0$ with such kernels) such that $AB=BF(A)$. In fact, let $n$ be a positive integer and consider $b(t)=\sum\limits_{j=1}^n b_j t^{-j}$, where $t\in [M_1,M_2]$, $b_j\in\mathbb{R}$ for $j=1,\ldots,n$ and $b_n\not=0$. We set $k_1(t,s)$ as defined by \eqref{KOneRemCorAIntOpBIndentOP}.
Then we have
\begin{align*}
& b(t)F_2(k(t,s))=b(t)(\delta_1 k(t,s)+\delta_2 k_1(t,s) )=(a_0 \delta_1 +\delta_2 \nu_0)\sum_{j=1}^n b_j t^{-j} \\
               &+(a_1\delta_1+\delta_2 \nu_1) \sum_{j=1}^n b_j t^{-j+1}+
              ( c_1\delta_1+\delta_2 \nu_2) \sum_{j=1}^n b_j t^{-j}s
              + \nu_3 \delta_2 \sum_{j=1}^n b_j t^{-j+1}s\\
              &
              = (a_1\delta_1+\delta_2\nu_1)b_1+\nu_3\delta_2b_1s+\sum_{j=1}^{n-1} ((a_0 \delta_1 +\delta_2 \nu_0)b_j +(a_1\delta_1+\delta_2 \nu_1)b_{j+1})t^{-j}\\
              &\resizebox{0.97\hsize}{!}{$\displaystyle +(a_0 \delta_1 +\delta_2 \nu_0)b_{n}t^{-n}+ \sum_{j=1}^{n-1} (( c_1\delta_1+\delta_2 \nu_2)b_j +\nu_3 \delta_2b_{j+1})t^{-j}s+( c_1\delta_1+\delta_2 \nu_2)b_nt^{-n}s $}
\\[0,2cm]
& k(t,s)b(s)=\,a_0\sum\limits_{j=1}^n b_j s^{-j}+a_1 \sum\limits_{j=1}^n b_j s^{-j} t +c_1 \sum\limits_{j=1}^n b_j s^{-j+1}\\
=&\,c_1b_1+\sum\limits_{j=1}^{n-1} (a_0b_j +c_1b_{j+1})s^{-j}+\sum\limits_{j=1}^{n} a_1b_j s^{-j} t+a_0b_ns^{-n}.
\end{align*}
 Thus we have
 $ k(t,s)b(s)=b(t)F_2(k(t,s))$ for almost every $(t,s)\in [M_1,M_2]\times[\alpha,\beta]$ if and only if
 \begin{eqnarray} 
 \nonumber
 c_1b_1&=& a_1\delta_1 b_1 + \delta_2 \nu_1 b_1, \\ \nonumber
 0&=&\delta_2 \nu_3 b_1,\\ \nonumber
 0&=&  (a_0\delta_1+\delta_2\nu_0) b_{j}+(\delta_1 a_1  +\delta_2 \nu_1) b_{j+1},\quad 1\le j\le n-1,\\
 \label{EqSevenRemLaurentCorIntOPMultOpPolinomials}
 a_0 b_j+c_1 b_{j+1}&=&0,\quad 1\le j \le n-1, \\
 \nonumber
 0&=&  c_1\delta_1b_j+\delta_2\nu_2 b_{j}+\delta_2 \nu_3 b_{j+1},\quad 1\le j\le n-1,\\
 \label{EqElevenRemCorLaurentIntOpMultOpPolinomials}
  a_1 b_j&=&0, \quad 1\le j\le n, \\
 \nonumber
 0&=& a_0\delta_1b_n+\delta_2 \nu_0 b_n,\\ \label{EqazerobnRemLaurentPolyIntOpMultop}
 0&=& a_0b_n,
\\ \nonumber
 0 &=& c_1\delta_1 b_n + \delta_2 \nu_3 b_n.
 \end{eqnarray}
Since $b_n\not=0$ then
from the equation \eqref{EqazerobnRemLaurentPolyIntOpMultop} we have $a_0b_n=0$ and thus $a_0=0$. From the equation \eqref{EqSevenRemLaurentCorIntOPMultOpPolinomials} for $j=n-1$ we get $c_1b_n=0$ and thus $c_1=0$. Finally from the equation \eqref{EqElevenRemCorLaurentIntOpMultOpPolinomials} we have
$0= a_1b_j $ for $j=n$ and thus $a_1=0$. This implies that $k(t,s)\equiv 0$, that is, $A= 0$. So $b_n\not= 0$ implies
$A= 0$. }
\end{rem:DjinjaTumwSilv1}

\begin{cor:DjinjaTumwSilv1}
Let $A:L_p(\mathbb{R})\to L_p(\mathbb{R})$,  $B:L_p(\mathbb{R})\to L_p(\mathbb{R})$, $1<p<\infty$,  be defined as follows, for almost all $t$,
\begin{gather*}
  (Ax)(t)= \int\limits_{\alpha}^{\beta}k(t,s) x(s)ds,\quad (Bx)(t)= b(t) x(t), \quad \alpha,\beta\in \mathbb{R}, \alpha<\beta,
\end{gather*}
where $k(t,s):\mathbb{R}\times [\alpha,\beta]\to \mathbb{R}$ is a measurable function, and $b\in L_\infty(\mathbb{R})$ is a nonzero function such that the set
$
 \supp  b(t)\cap [\alpha,\beta]
$
has measure zero.
Consider a polynomial $F(z)=\sum\limits_{j=0}^{n}\delta_j z^j$, where $\delta_0,\ldots,\delta_n\in \mathbb{R}$. Let
\begin{gather*}\nonumber
  k_0(t,s)=k(t,s), \quad
 \quad k_m(t,s)=\int\limits_{\alpha}^{\beta} k(t,\tau)k_{m-1}(\tau,s)d\tau,\quad m=1,\ldots,n,
  \\
 F_0(k(t,s))=0,\quad  F_n(k(t,s))=\sum_{j=1}^{n} \delta_j k_{j-1}(t,s),\quad n\geq 1.
\end{gather*}
Then $ AB=BF(A)$
if and only if $\delta_0=0$ and the set
\begin{gather*}
  \supp b(t)\,\cap\, \supp F_n(k(t,s))
\end{gather*}
has measure zero in $\mathbb{R}\times [\alpha,\beta]$.
\end{cor:DjinjaTumwSilv1}

\begin{proof}
Suppose that  the set
$
 \supp  b\cap [\alpha,\beta]
$
has measure zero.
By Proposition \ref{propIntOpLpAintBidentity} we have $AB=BF(A)$ if and only if condition  \eqref{CondABeqBFADeltaZeroIntOpIdentity} holds, that is,
\begin{gather*}
\forall \ x\in L_p(\mathbb{R}): \quad  b(t)\delta_0 x(t)+
   b(t)\int\limits_{\alpha}^{\beta} F_n(k(t,s))x(s)ds = \int\limits_{\alpha}^{\beta}  k(t,s)b(s) x(s)ds,
\end{gather*}
almost everywhere. By taking  $x(\cdot)=I_{[M_1,M_2]}(\cdot)b(\cdot)$, where $M_1,M_2\in\mathbb{R}$, $M_1<M_2$, $[M_1,M_2]\supset [\alpha,\beta]$, $\tilde{m}([M_1,M_2]\setminus [\alpha,\beta])>0$, $I_E(\cdot)$ is the indicator function of the set $E$, the condition \eqref{CondABeqBFADeltaZeroIntOpIdentity} reduces to
\begin{equation*}
I_{[M_1,M_2]}(\cdot) b(\cdot)^2\delta_0 =0.
\end{equation*}
Since $b$ has support with positive measure (otherwise $B\equiv 0$), then $\delta_0=0$. By using this, condition
\eqref{CondABeqBFADeltaZeroIntOpIdentity} reduces to the following
\begin{gather*}
\forall \ x\in L_p(\mathbb{R}): \quad
   b(t)\int\limits_{\alpha}^{\beta} F_n(k(t,s))x(s)ds = \int\limits_{\alpha}^{\beta}  k(t,s)b(s) x(s)ds.
\end{gather*}
By hypothesis the right hand side is equal zero. Then condition \eqref{CondABeqBFADeltaZeroIntOpIdentity} reduces to
\begin{gather*}
\forall \ x\in L_p(\mathbb{R}): \quad
   b(t)\int\limits_{\alpha}^{\beta} F_n(k(t,s))x(s)ds =0.
\end{gather*}
This is equivalent to
\begin{gather} \label{ProofCorSuppBOutInterval}
b(t)F_n(k(t,s))=0 \quad \mbox{ for almost every }\ s\in[\alpha,\beta].
\end{gather}
By applying a similar argument
as in the proof of Proposition \ref{propIntOpLpAintBidentity}
we conclude that condition \eqref{ProofCorSuppBOutInterval} is equivalent to that the set
\begin{gather*}
  \supp b(t)\,\cap\, \supp F_n(k(t,s))
\end{gather*}
has measure zero in $\mathbb{R}\times [\alpha,\beta]$.
\qed \end{proof}

\begin{cor:DjinjaTumwSilv1}\label{corIntOpLpAintBidentity}
Let $A:L_p(\mathbb{R})\to L_p(\mathbb{R})$,  $B:L_p(\mathbb{R})\to L_p(\mathbb{R})$, $1\leq p \leq\infty$, be defined as follows, for almost all $t$,
\begin{gather*}
  (Ax)(t)= \int\limits_{\alpha}^{\beta}a(t)c(s) x(s)ds,\quad (Bx)(t)= b(t) x(t), \quad \alpha,\beta\in \mathbb{R}, \alpha<\beta,
\end{gather*}
where $a:\mathbb{R}\to\mathbb{R}$, $c:[\alpha,\beta]\to \mathbb{R}$, $b:\mathbb{R}\to \mathbb{R}$ are measurable functions.
Consider a polynomial defined by $F(z)=\sum\limits_{j=1}^{n}\delta_j z^j$, where $\delta_1,\ldots,\delta_n\in \mathbb{R}$. Let
$
\mu=\int\limits_{\alpha}^{\beta} a(s)c(s)ds.
$
Then, $ AB=BF(A)$ if and only if the set
\begin{gather*}
   \supp (a(t)c(s))\cap \supp \left( b(t) \sum\limits_{j=1}^{n} \delta_j \mu^{j-1}-b(s)\right)
 \end{gather*}
has measure zero in $\mathbb{R}\times[\alpha,\beta]$.
\end{cor:DjinjaTumwSilv1}

\begin{proof}
 We set $k(t,s)=a(t)c(s)$,   so we  have
 \begin{gather*}
   k_0(t,s)=k(t,s)=a(t)c(s),\quad
   \\
   k_m(t,s)=\int\limits_{\alpha}^{\beta} k(t,\tau)k_{m-1}(\tau,s)d\tau =a(t)c(s)\left(\int\limits_{\alpha}^{\beta} a(s)c(s)ds\right)^{m}, \ m=1,\ldots,n\\
   F_n(k(t,s))=\sum_{j=1}^{n} \delta_j k_{j-1}(t,s)=
   \sum_{j=1}^{n} \delta_j a(t)c(s)\left(\int\limits_{\alpha}^{\beta} a(s)c(s)ds\right)^{j-1},\  \quad n\geq 1.
 \end{gather*}
 By applying Proposition \ref{propIntOpLpAintBidentity} we have  $AB=BF(A)$ if and only if
 \begin{gather*}
 b(t) \sum\limits_{j=1}^{n} \delta_j a(t)c(s)\left(\int\limits_{\alpha}^{\beta} a(s)c(s)ds\right)^{j-1}=a(t)c(s)b(s)  \
 \Longleftrightarrow \\
 a(t)c(s)\left( b(t) \sum\limits_{j=1}^{n} \delta_j \left(\int\limits_{\alpha}^{\beta} a(s)c(s)ds\right)^{j-1}-b(s)\right)=0
 \end{gather*}
 for almost every $(t,s)$ in $\mathbb{R}\times [\alpha,\beta]$.
 The last condition is equivalent to the set
 \begin{gather*}
   \supp (a(t)c(s))\cap \supp \left( b(t) \sum\limits_{j=1}^{n} \delta_j \left(\int\limits_{\alpha}^{\beta} a(s)c(s)ds\right)^{j-1}-b(s)\right)
 \end{gather*}
 has measure zero in $\mathbb{R}\times[\alpha,\beta]$. We complete the proof by noticing that the corresponding  set can be written as
 \begin{gather*}
   \supp (a(t)c(s))\cap \supp \left( b(t) \sum\limits_{j=1}^{n} \delta_j \mu^{j-1}-b(s)\right),
 \end{gather*}
 where $\displaystyle \mu=\int\limits_{\alpha}^{\beta} a(s)c(s)ds$.
  \qed \end{proof}


\begin{ex:DjinjaTumwSilv1} \label{ExampleIntOpMultOpCondOnPolinomial} {\rm
Let $A:L_p(\mathbb{R})\to L_p(\mathbb{R})$, $B:L_p(\mathbb{R})\to L_p(\mathbb{R})$, $1<p<\infty$ be  defined as follows, for almost all $t$,
$$
  (Ax)(t)= \int\limits_{0}^{2} a(t)c(s) x(s)ds,\ (Bx)(t)= b(t)x(t),
$$
where $a(t)=2tI_{[0,2]}(t) $, $c(s)=I_{[0,1]}(s)$, $b(t)=I_{[1,2]}(t) t^2 $.
 Since kernel has compact support, we can apply \cite[Theorem 3.4.10]{HutsonPymDST} and, we conclude that  operators $A$ is well defined and bounded. Since function $b$ has $4$ as an upper bound then $\|B\|_{L_p} \leq 4$. Hence operator $B$ is well defined and bounded.
 Consider a polynomial defined by $F(z)=\sum\limits_{j=1}^{n}\delta_j z^j$, where $\delta_1,\ldots,\delta_n\in \mathbb{R}$. Then, the above operators satisfy the relation
$ AB=BF(A)$ if and only if $\sum\limits_{j=1}^n \delta_j=0$.  In fact, by applying Corollary \ref{corIntOpLpAintBidentity} we have
 \begin{gather*}
    \mu=\int\limits_{0}^2 a(s)c(s)ds=1.
 \end{gather*}
Hence,
$\supp \{b(t)\cdot 0-b(s) \}=\mathbb{R}\times [1,2].$
Moreover,
$\supp a(t)c(s)=[0,2]\times [0,1].$
The set
$\supp (a(t)c(s))\cap \supp \left( -b(s)\right),$
has measure zero in $\mathbb{R}\times [0,2]$.
}
 \end{ex:DjinjaTumwSilv1}

 \begin{ex:DjinjaTumwSilv1}{\rm
Let $A:L_p(\mathbb{R})\to L_p(\mathbb{R})$,  $B:L_p(\mathbb{R})\to L_p(\mathbb{R})$, $1<p<\infty$ be  defined as follows, for almost all $t$,
\begin{gather*}
  (Ax)(t)= \int\limits_{0}^{2} a(t)c(s) x(s)ds,\quad (Bx)(t)= b(t)x(t),
\end{gather*}
where $a(t)=I_{[0,2]}(t)\sin(\pi t) $, $c(s)=I_{[0,1]}(s)$, $b(t)=I_{[1,2]}(t) t^2 $.
  Since $a\in L_p(\mathbb{R})$ and $c\in L_q[0,2]$, $1<q<\infty$, $\frac{1}{p}+\frac{1}{q}=1$, by applying H\"older inequality we have that  operator $A$ is well defined and bounded. The function $b\in L_{\infty}$, so $B$ is well defined and bounded because $\|B\|_{L_p}\leq \|b\|_{L_\infty}$ we conclude that operator
 $B$ is well defined and bounded.
 Consider a polynomial defined by $F(z)=\delta z^d$, where $\delta\not=0$ is a real constant and $d$ is a positive integer $d\ge 2$. Then, the above operators satisfy the relation
$
  AB=\delta BA^d.
$
  In fact, by applying Corollary \ref{corIntOpLpAintBidentity} we have
$
    \mu=\int\limits_{0}^2 a(s)c(s)ds=0.
$
Hence,
$   \supp \{b(t)\cdot 0-b(s) \}=\mathbb{R}\times [1,2].$
Moreover,
$\supp a(t)c(s)=[0,2]\times [0,1].$
The set
$\supp (a(t)c(s))\cap \supp \left( -b(s)\right),$
has measure zero in $\mathbb{R}\times [0,2]$.
}
\end{ex:DjinjaTumwSilv1}

\begin{ex:DjinjaTumwSilv1}{\rm
Let $A:L_p(\mathbb{R})\to L_p(\mathbb{R})$,  $B:L_p(\mathbb{R})\to L_p(\mathbb{R})$, $1<p<\infty$, be defined as follows, for almost all $t$,
\begin{gather*}
  (Ax)(t)= \int\limits_{\alpha}^{\beta} I_{[\alpha,\beta]}(t) x(s)ds,\quad (Bx)(t)= I_{[\alpha,\beta]}(t) x(t), \quad \alpha,\beta\in\mathbb{R}, \alpha<\beta.
\end{gather*}
  Since kernel has compact support, we can apply \cite[Theorem 3.4.10]{HutsonPymDST} and, we conclude that  operator $A$ is well defined and bounded. Since  $\|B\|_{L_p}\leq 1$ then operator  $B$ is well defined and bounded.
 Consider a polynomial defined by $F(z)=\sum\limits_{j=1}^{n}\delta_j z^j$, where $\delta_1,\ldots,\delta_n\in \mathbb{R}$. Then, the above operators satisfy the relation
$AB=BF(A)$ if and only if
$\sum\limits_{j=1}^{n} \delta_j (\beta-\alpha)^{j-1}=1.$
Indeed, if $a(t)=b(t)=I_{[\alpha,\beta]}(t)$, $c(s)=1$ and
\begin{gather*}
  \lambda=\sum\limits_{j=1}^{n} \delta_j \left(\int\limits_{\alpha}^{\beta}a(s)c(s)ds\right)^{j-1}=\sum\limits_{j=1}^{n} \delta_j (\beta-\alpha)^{j-1},
\end{gather*}
then from Corollary \ref{corIntOpLpAintBidentity} we have the following:
\begin{itemize}
  \item If $\lambda\not=0$, $\lambda\not=1$,
\begin{gather*}
  \supp (b(t)\lambda-b(s)))=\left\{(t,s)\in\mathbb{R}\times [\alpha,\beta]:\ \lambda I_{[\alpha,\beta]}(t)\not=1\right\}=\mathbb{R}\times[\alpha,\beta],\\ \\
  \supp  a(t)c(s)=\{(t,s)\in\mathbb{R}\times [\alpha,\beta]:\ I_{[\alpha,\beta]}(t)\not=0\}=[\alpha,\beta]\times [\alpha,\beta].
\end{gather*}
The set
$ \supp (\lambda b(t)-b(s))\cap \supp (a(t)c(s))=[\alpha,\beta]\times [\alpha,\beta]
$
has  positive measure.

\item If $\lambda=1$,
\begin{eqnarray*}
  \supp (b(t)-b(s))&=&\left\{(t,s)\in\mathbb{R}\times [\alpha,\beta]:\ I_{[\alpha,\beta]}(t)\not=1\right\}\\
  &=&(\mathbb{R}\setminus[\alpha,\beta])\times [\alpha,\beta].
\end{eqnarray*}
The set
$
  \supp ( b(t)-b(s))\cap \supp (a(t)c(s))
$
has measure zero in $\mathbb{R}\times [\alpha,\beta]$.
\item If $\lambda=0$,
\begin{eqnarray*}
  \supp (\lambda b(t)-b(s))&=&\supp b(s)=\left\{(t,s)\in\mathbb{R}\times[\alpha,\beta]:\ I_{[\alpha,\beta]}(s)\not=0\right\}\\
  &=& \left\{(t,s)\in\mathbb{R}\times[\alpha,\beta]:\ \alpha\le s\le \beta \right\}.
\end{eqnarray*}
The set
$
  \supp  b(s)\cap \supp (a(t)c(s))=[\alpha,\beta]\times [\alpha,\beta]
$
has measure $(\beta-\alpha)^2$.
\end{itemize}
 The conditions in the Corollary \ref{corIntOpLpAintBidentity} are fulfilled only in the second case, that is, when $\lambda=1$.
}
\end{ex:DjinjaTumwSilv1}

\subsubsection{Representations when \texorpdfstring{$A$}{A} is multiplication operator and \texorpdfstring{$B$}{B} is integral operator}

\begin{prop:DjinjaTumwSilv1}\label{propIntOpBLpIdentityA}
Let $A:L_p(\mathbb{R})\to L_p(\mathbb{R})$,  $B:L_p(\mathbb{R})\to L_p(\mathbb{R})$, $1<p<\infty$ be
defined as follows, for almost all $t$,
\begin{gather*}
  (Ax)(t)=a(t) x(t),\quad  (Bx)(t)=\int\limits_{\alpha}^{\beta}k(t,s) x(s)ds, \quad \alpha,\beta\in \mathbb{R}, \alpha<\beta,
\end{gather*}
where $a:\mathbb{R}\to \mathbb{R}$,  $k: \mathbb{R}\times [\alpha,\beta]\to \mathbb{R}$ are measurable functions. Consider a polynomial defined by $F(z)=\sum\limits_{j=0}^{n}\delta_j z^j$, where $\delta_0,\ldots,\delta_n\in \mathbb{R}$. Then
\begin{gather*}
  AB=BF(A)
\end{gather*}
if and only if  the set
   \begin{gather*}
    \supp (a(t)-F(a(s)))\cap \supp k(t,s)
      \end{gather*}
  has measure zero in $\mathbb{R}\times [\alpha,\beta]$.
\end{prop:DjinjaTumwSilv1}

\begin{proof}
  We have for almost every $t\in\mathbb{R}$
  \begin{gather*}
    (ABx)(t)= \int\limits_{\alpha}^{\beta} a(t) k(t,s) x(s)ds \\ 
    (A^nx)(t)=[a(t)]^nx(t) \\
    (F(A)x)(t)= \sum_{i=0}^{n} \delta_i (A^ix)(t)=\left(\sum_{i=0}^{n} \delta_i[a(t)]^i \right) x(t)=F(a(t))x(t)\\ 
    (BF(A)x)(t)=\int\limits_{\alpha}^{\beta}k(t,s))F(a(s)) x(s)ds.
  \end{gather*}

    Then we have $ABx=BF(A)x$ if and only if
  \begin{gather}\label{IntEqOpProofPropInOpIdentOp}
    \int\limits_{\alpha}^{\beta}a(t) k(t,s) x(s)ds=\int\limits_{\alpha}^{\beta} k(t,s)F(a(s)) x(s)ds.
  \end{gather}
  almost everywhere. By using Lemma \ref{lemEqIntForAllLpFunct} and by applying the same argument
  as in the final steps on the proof of Proposition \ref{propIntOpLpAintBidentity},  the condition \eqref{IntEqOpProofPropInOpIdentOp} is equivalent to
      \begin{gather*}
         a(t)k(t,s)=k(t,s)F(a(s)) \Longleftrightarrow k(t,s)(a(t)-F(a(s)))=0
          \end{gather*}
        for  almost every $(t,s)$ in $\mathbb{R}\times [\alpha,\beta]$.

   Since the variables $t$ and $s$ are independent, this is true if and only if the set
           \begin{gather*}
    \supp (a(t)-F(a(s)))\cap \supp k(t,s)
      \end{gather*}
  has measure zero in $\mathbb{R}\times [\alpha,\beta]$.
  \qed \end{proof}

\begin{ex:DjinjaTumwSilv1}{\rm
 Let $A:L_p(\mathbb{R})\to L_p(\mathbb{R})$,  $B:L_p(\mathbb{R})\to L_p(\mathbb{R})$, $1<p<\infty$
be defined as follows, for almost all $t$,
\begin{gather*}
(Ax)(t)=I_{[\alpha,\beta]}(t) x(t),\quad  (Bx)(t)=\int\limits_{\alpha}^{\beta} I_{[\alpha,\beta]^2}(t,s)x(s)ds, \quad \alpha,\beta\in \mathbb{R}, \alpha <\beta
\end{gather*}
By using properties of norm and \cite[Theorem 3.4.10]{HutsonPymDST}, respectively, for operators $A$ and $B$, we conclude that operators
 $A$ and $B$ are well defined and bounded. 
 For a monomial defined by $F(z)=z^n$, $n=1,2,\ldots$, the above operators satisfy the relation
$
  AB=BF(A).
$
In fact, by setting $a(t)=I_{[\alpha,\beta]}(t)$, $k(t,s)=I_{[\alpha,\beta]^2}(t,s)$ we have
\begin{gather*}
  \supp (a(t)-F(a(s)))=\left\{(t,s)\in\mathbb{R}\times [\alpha,\beta]:\ I_{[\alpha,\beta]}(t)\not=1\right\}=(\mathbb{R}\setminus [\alpha,\beta])\times [\alpha,\beta],\\ \\
  \supp  k(t,s)=[\alpha,\beta]\times [\alpha,\beta].
\end{gather*}
The set
$\supp (a(t)-F(a(s)))\cap \supp (k(t,s))$
has measure zero in $\mathbb{R}\times[\alpha,\beta]$. So the result follows from  Proposition \ref{propIntOpBLpIdentityA}.
}\end{ex:DjinjaTumwSilv1}

\begin{ex:DjinjaTumwSilv1}{\rm
 Let $A:L_p(\mathbb{R})\to L_p(\mathbb{R})$,  $B:L_p(\mathbb{R})\to L_p(\mathbb{R})$, $1<p<\infty$
defined as follows, for almost all $t$,
\begin{equation*}
  (Ax)(t)=(\gamma_1I_{[0,1/2)}(t)-\gamma_2I_{[1/2,1]}(t)) x(t)                ,\quad  (Bx)(t)=\int\limits_{0}^{1} k(t,s) x(s)ds
\end{equation*}
$k:\mathbb{R}\times [0,1]\to \mathbb{R}$ is a Lebesgue measurable function such that $B$ is well defined. The operator $A$ is well defined and bounded.
Consider a polynomial defined by $F(z)=\delta_0+\delta_1 z$,  where $\delta_0$, $\delta_1$, $\gamma_1$, $\gamma_2$ are constants such that
\begin{gather*}
|\delta_0|+|\delta_1|+|\gamma_1|+|\gamma_2|\not=0.
\end{gather*}
 If $k(\cdot,\cdot)$ is a measurable function such that one of the following is fulfilled:
\begin{enumerate}
  \item $\delta_0=-\delta_1\gamma_1$ and $\supp \ k(t,s)\subseteq (\mathbb{R}\setminus [0,1])\times [0,1/2]$;

    \item $\delta_0=\delta_1\gamma_2$ and $\supp \ k(t,s)\subseteq (\mathbb{R}\setminus [0,1])\times [1/2,1]$;



          \item $\delta_0+\delta_1 \gamma_1-\gamma_1=0$ and $\supp  k(t,s)\subseteq [0,1/2]\times [0,1/2]$;

            \item $\delta_0+\delta_1 \gamma_1+\gamma_2=0$ and $\supp k(t,s)\subseteq [1/2,1]\times [0,1/2]$;

              \item $\delta_0-\delta_1 \gamma_2-\gamma_1=0$ and $\supp  k(t,s)\subseteq [0,1/2]\times [1/2,1]$;

                \item $\delta_0-\delta_1 \gamma_2+\gamma_2=0$ and $\supp  k(t,s)\subseteq [1/2,1]\times [1/2,1]$,
\end{enumerate}
then the above operators satisfy the relation $AB=BF(A)$.

In fact, putting $a(t)=\gamma_1 I_{[0,1/2)}(t)-\gamma_2 I_{[1/2,1]}(t)$ we have
\begin{gather*}
 [a(t)-F(a(s))]=\left\{\begin{array}{cccc}
                                0, & \mbox{ if } \delta_0=-\delta_1\gamma_1, & t\not\in[0,1], & s\in[0,1/2) \\
                                0, & \mbox{ if } \delta_0=\delta_1\gamma_2, & t\not\in[0,1], & s\in[1/2,1] \\
                                0, & \mbox{ if } \delta_0+\delta_1\gamma_1-\gamma_1=0, & t\in[0,1/2), & s\in[0,1/2) \\
                                 0, & \mbox{ if } \delta_0+\delta_1\gamma_1+\gamma_2=0, & t\in[1/2,1), & s\in[0,1/2] \\
                                 0, & \mbox{ if } \delta_0-\delta_1\gamma_2-\gamma_1=0, & t\in[0,1/2], & s\in[1/2,1] \\
                                0, & \mbox{ if } \delta_0-\delta_1\gamma_2+\gamma_2=0, & t\in[1/2,1], & s\in[1/2,1] \\
                                \gamma_3, & \mbox{ otherwise }  &  &  \\
                                \end{array}\right.
\end{gather*}
where $\gamma_3$ can be different from zero depending on the constants involved. Thus, in each condition we can choose
$k(t,s)=I_{S}(t,s)$, where $S=\{(t,s)\in\mathbb{R}\times [0,1]:\ a(t)-F(a(s))=0\}$ and with a positive measure. Or for instance we can take:
\begin{enumerate}
  \item  $k(t,s)=I_{[2,3]\times [0,1/2]}(t,s)$ if $\delta_0=-\delta_1\gamma_1$;
  \item  $k(t,s)=I_{[2,3]\times [1/2,1]}(t,s)$ if $\delta_0=\delta_1\gamma_2$;
  \item  $k(t,s)=I_{[0,1/3]\times [1/3,1/2]}(t,s)$ if $\delta_0+\delta_1\gamma_1-\gamma_1=0$;
  \item  $k(t,s)=I_{[2/3,1/2]\times [0,1/2]}(t,s)$ if $\delta_0+\delta_1\gamma_1+\gamma_2=0$;
  \item  $k(t,s)=I_{[0,1/3]\times [2/3,1]}(t,s)$ if $\delta_0-\delta_1\gamma_2-\gamma_1=0$;
  \item  $k(t,s)=I_{[2/3,1]\times [2/3,1]}(t,s)$ if $\delta_0-\delta_1\gamma_2+\gamma_2$.
\end{enumerate}

According to the definition, in all above cases the set
\begin{equation*}
  \supp (a(t)-F(a(s)))\cap \supp\ (k(t,s))
\end{equation*}
has measure zero in $\mathbb{R}\times [0,1]$. So the result follows from  Proposition \ref{propIntOpBLpIdentityA}.
}\end{ex:DjinjaTumwSilv1}

\begin{cor:DjinjaTumwSilv1}\label{CorIntOpAidentBIntOpSplittedK}
Let $A:L_p(\mathbb{R})\to L_p(\mathbb{R})$,  $B:L_p(\mathbb{R})\to L_p(\mathbb{R})$, $1<p<\infty$
defined as follows, for almost all $t$,
\begin{gather*}
  (Ax)(t)=a(t) x(t),\quad  (Bx)(t)=\int\limits_{\alpha}^{\beta}b(t)c(s) x(s)ds, \quad \alpha,\beta\in \mathbb{R}, \alpha<\beta,
\end{gather*}
where $a:\mathbb{R}\to\mathbb{R}$, $b:\mathbb{R}\to\mathbb{R}$, $c:[\alpha,\beta]\to\mathbb{R}$ are measurable functions. For a polynomial defined by $F(z)=\sum\limits_{j=0}^{n}\delta_j z^j$, where $\delta_0,\ldots,\delta_n\in \mathbb{R}$, we have
\begin{gather*}
  AB=BF(A)
\end{gather*}
if and only if  the set
   \begin{equation*}
    \supp (a(t)-F(a(s)))\cap \supp (b(t)c(s))
      \end{equation*}
  has measure zero in $\mathbb{R}\times [\alpha,\beta]$.
\end{cor:DjinjaTumwSilv1}

\begin{proof}
This follows by Proposition \ref{propIntOpBLpIdentityA}.
\qed \end{proof}

\begin{ex:DjinjaTumwSilv1}{\rm
 Let $A:L_p(\mathbb{R})\to L_p(\mathbb{R})$,  $B:L_p(\mathbb{R})\to L_p(\mathbb{R})$, $1<p<\infty$
be defined as follows, for almost all $t$,
\begin{equation*}
  (Ax)(t)=a(t) x(t),\quad  (Bx)(t)=\int\limits_{\alpha}^{\beta} b(t)c(s) x(s)ds,
 \quad  \alpha,\beta\in \mathbb{R}, \alpha<\beta,
\end{equation*}
where $a(t)=-1 +I_{[\alpha,\beta]}(t)$, $b(t)=I_{[\alpha-2,\alpha-1]}(t)$, $c(s)=1$.
We have that $a\in L_{\infty}(\mathbb{R})$ and so $\| A\|_{L_p}\leq \|a\|_{L_\infty}$. Therefore, $A$ is well defined and bounded. Since kernel has compact support in $\mathbb{R}\times [\alpha,\beta]$, we can apply \cite[Theorem 3.4.10]{HutsonPymDST} and, we conclude that  operators $B$ is well defined and bounded.
Consider a polynomial defined by $F(z)=-1+\delta_1 z$, where $\delta_1$ is a real constant. Then the above operators satisfy the relation
$
  AB=BF(A).
$
In fact, for $(t,s)\in\mathbb{R}\times [\alpha,\beta]$  we have
\begin{eqnarray*}
  F(a(s))-a(t)&=&-\delta_1+ \delta_1 I_{[\alpha,\beta]}(s)-I_{[\alpha,\beta]}(t)=
  -I_{[\alpha,\beta]}(t).
\end{eqnarray*}
Therefore, we have
\begin{gather*}
  \supp (a(t)-F(a(s)))=[\alpha,\beta]\times [\alpha,\beta],\\
  \supp  b(t)c(s)=\supp \ I_{[\alpha-2,\alpha-1]}(t)I_{[\alpha,\beta]}(s)=[\alpha-2,\alpha-1]\times [\alpha,\beta].
\end{gather*}
The set
$
  \supp (a(t)-F(a(s)))\cap \supp (I_{[\alpha-2,\alpha-1]}(t)I_{[\alpha,\beta]}(s))
$
has measure zero. So the result follows from  Corollary \ref{CorIntOpAidentBIntOpSplittedK}.
}\end{ex:DjinjaTumwSilv1}

\begin{ex:DjinjaTumwSilv1}{\rm
 Let $A:L_p(\mathbb{R})\to L_p(\mathbb{R})$,  $B:L_p(\mathbb{R})\to L_p(\mathbb{R})$, $1<p<\infty$
be defined as follows, for almost all $t$,
\begin{gather*}
  (Ax)(t)=a(t) x(t),\quad  (Bx)(t)=\int\limits_{\alpha}^{\beta} b(t)c(s) x(s)ds, \quad \alpha,\beta\in \mathbb{R}, \alpha<\beta,
\end{gather*}
where $a(t)=\gamma_0+I_{\left[\alpha, \frac{\alpha+\beta}{2}\right]}(t)t^2$, $\gamma_0$ is a real number, $b(t)=(1+t^2)I_{[\beta+1,\beta+2]}(t) $, $c(s)=I_{\left[\frac{\alpha+\beta}{2},\beta\right]}(s)(1+s^4)$. Consider a polynomial defined by $F(z)=\delta_0+\delta_1 z$, where $\delta_0,\delta_1$ are real constants and $\delta_1\not=0$. If $\delta_0=\gamma_0-\delta_1 \gamma_0$ then the above operators satisfy the relation
\begin{gather*}
  AB-\delta_1 BA=\delta_0 B.
\end{gather*}
In fact, $A$ is well defined, bounded since $a\in L_\infty$ and this implies $\|A\|_{L_p}\leq \|a\|_{L_\infty}$. Operator $B$ is well defined, bounded since $k(t,s)=b(t)c(s)$, $(t,s)\in \mathbb{R}\times [\alpha,\beta]$ has compact support and satisfies conditions of \cite[Theorem 3.4.10]{HutsonPymDST}. If $\delta_0=\gamma_0-\delta_1 \gamma_0$ then we have
\begin{eqnarray*}
  F(a(s))-a(t)&=& \delta_0+\gamma_0\delta_1+ \delta_1 I_{\left[\alpha,\frac{\alpha+\beta}{2}\right]}(s)s^2-\gamma_0-I_{\left[\alpha,\frac{\alpha+\beta}{2}\right]}(t)t^2\\
     &=& \delta_1 I_{\left[\alpha, \frac{\alpha+\beta}{2}\right]}(s)s^2-I_{\left[\alpha,\frac{\alpha+\beta}{2}\right]}(t)t^2.
  \end{eqnarray*}
Then we have
$$
  \supp [a(t)-F(a(s))]=
  \left( \mathbb{R}\times \left[\alpha,\frac{\alpha+\beta}{2}\right]\cup \left[\alpha,\frac{\alpha+\beta}{2}\right]\times\left[\frac{\alpha+\beta}{2},\beta\right]\right)\setminus W,
$$
  where $W\subseteq \mathbb{R}\times [\alpha,\beta]$ is a set with Lebesgue measure zero, and
   \begin{eqnarray*}
  \supp \ b(t)c(s)&=&\supp \ (1+t^2)I_{[\beta+1,\beta+2]}(t) I_{\left[\frac{\alpha+\beta}{2},\beta\right]}(s)(1+s^4)\\
  &=&[\beta+1,\beta+2]\times \left[\frac{\alpha+\beta}{2},\beta\right].
\end{eqnarray*}
The set
$
  \supp (a(t)-F(a(s)))\cap \supp (b(t)c(s))
$
has measure zero. So the result follows from  Corollary \ref{CorIntOpAidentBIntOpSplittedK}.
}
\end{ex:DjinjaTumwSilv1}


\acknowledgement{
This work was supported by the Swedish International Development Cooperation
Agency (Sida) bilateral program with  Mozambique. Domingos Djinja is grateful to the research environment Mathematics and Applied Mathematics (MAM), Division of Mathematics and Physics, School of Education, Culture and Communication, M\"alardalen University for excellent environment for research in Mathematics.
Partial support from Swedish Royal Academy of Sciences is also gratefully acknowledged.}

\end{document}